\input amstex \documentstyle{amsppt} \nologo
\loadbold
\let\bk\boldkey

\hsize=5.75truein
\vsize=8.75truein 
\hcorrection{.25truein}
\loadeusm \let\scr\eusm
\loadeurm 
\font\Smc=cmcsc10 scaled \magstep1
\define\cind{\text{\it c-\/\rm Ind}} 
\define\Ind{\text{\rm Ind}}
\define\Aut#1#2{\text{\rm Aut}_{#1}(#2)}
\define\End#1#2{\text{\rm End}_{#1}(#2)}
\define\Hom#1#2#3{\text{\rm Hom}_{#1}({#2},{#3})} 
\define\GL#1#2{\roman{GL}_{#1}(#2)}
\define\M#1#2{\roman M_{#1}(#2)}
\define\Gal#1#2{\text{\rm Gal\hskip.5pt}(#1/#2)}
\define\Ao#1#2{\scr A_{#1}(#2)} 
\define\Go#1#2{\scr G_{#1}(#2)} 
\define\Gfin#1#2{\scr G^{\text{\rm fin}}_{#1}(#2)} 
\define\Afin#1#2{\scr A^{\text{\rm fin}}_{#1}(#2)}
\define\upr#1#2{{}^{#1\!}{#2}} 
\define\wt#1{\widetilde{#1}} 
\define\N#1#2{\text{\rm N}_{#1/#2}}
\let\ge\geqslant
\let\le\leqslant
\let\ups\upsilon
\let\vD\varDelta 
\let\vG\varGamma
\let\vL\varLambda 

\let\vP\varPi 
\let\vS\varSigma 
\let\vt\vartheta
\let\vT\varTheta
\define\wP#1{\widehat{\scr P}_{#1}} 
\define\Gg#1#2#3{\scr G_{#1}(#2;#3)} 
\define\Ggf#1#2#3{\scr G^{\text{\rm fin}}_{#1}(#2;#3)} 
\define\Aa#1#2#3{\scr A_{#1}(#2;#3)} 
\define\Aaf#1#2#3{\scr A^{\text{\rm fin}}_{#1}(#2;#3)} 
\define\Q{\,\overline{\!\bk Q}_\ell} 
\define\Z{\,\overline{\!\bk Z}_{\!\ell}} 
\define\F{\,\overline{\!\bk F}_{\!\ell}} 
\define\pl{\frak P_\ell} 
\let\Bbb\bk
\topmatter 
\title 
A congruence property of the local Langlands correspondence 
\endtitle 
\rightheadtext{Modular Langlands correspondence} 
\author 
Colin J. Bushnell and Guy Henniart 
\endauthor 
\leftheadtext{C.J. Bushnell and G. Henniart}
\affil 
King's College London and Universit\'e de Paris-Sud 
\endaffil 
\address 
King's College London, Department of Mathematics, Strand, London WC2R 2LS, UK. 
\endaddress
\email 
colin.bushnell\@kcl.ac.uk 
\endemail
\address 
Institut Universitaire de France et Universit\'e de Paris-Sud, Laboratoire de Math\'ematiques d'Orsay,
Orsay Cedex, F-91405; CNRS, Orsay cedex, F-91405. 
\endaddress 
\email 
Guy.Henniart\@math.u-psud.fr 
\endemail 
\date June 2011 \enddate 
\abstract 
Let $F$ be a non-Archimedean local field of residual characteristic $p$, and $\ell$ a prime number, $\ell \neq p$. We consider the Langlands correspondence, between irreducible, $n$-dimensional, smooth representations of the Weil group of $F$ and irreducible cuspidal representations of $\text{\rm GL}_n(F)$. We use an explicit description of the correspondence from an earlier paper, and otherwise entirely elementary methods, to show that it respects the relationship of congruence modulo $\ell$. The $\ell$-modular correspondence thereby becomes as effective as the complex one. 
\endabstract 
\endtopmatter 
\document \nopagenumbers
Let $\ell$ be a prime number and let $F$ be a non-Archimedean local field of residual characteristic $p\neq\ell$. We form the Weil group $\scr W_F$ of $F$, relative to some separable algebraic closure $\overline F/F$. Let $\bk E$ denote either the field $\Bbb C$ of complex numbers or an algebraic closure $\Q$ of the field $\Bbb Q_\ell$ of $\ell$-adic numbers. For an integer $n\ge 1$, let $\Go nF_{\bk E}$ denote the set of equivalence classes of irreducible, smooth $\bk E$-representations of $\scr W_F$ of dimension $n$. Likewise, let $\Ao nF_\bk E$ be the set of equivalence classes of irreducible, smooth, cuspidal $\bk E$-representations of the group $G = \GL nF$. 
\par
The Langlands correspondence gives a canonical bijection $\Go nF_\Bbb C \to \Ao nF_\Bbb C$ which we denote $\sigma \mapsto \upr L\sigma$. By choosing a field isomorphism $\alpha:\Bbb C\to \Q$, one may transport the Langlands correspondence to a bijection $\Go nF_{\Q} \to \Ao nF_{\Q}$. We continue to denote this $\sigma \mapsto \upr L\sigma$, eliding the fact that it may also depend on the choice of $\alpha$. 
\par
Let $\F$ be an algebraic closure of the field $\Bbb F_\ell$ of $\ell$ elements. If $\sigma\in \Go nF_{\Q}$ has determinant of finite order then $\sigma(\scr W_F)$ is finite and, following a standard technique from the representation theory of finite groups, one may attach to $\sigma$ an isomorphism class $[\sigma]^\roman{ss}_\ell$ of $n$-dimensional, smooth, semisimple $\F$-representations of $\scr W_F$. Likewise, if $\pi\in \Ao nF_{\Q}$ has central character of finite order, a variation on the same technique attaches to $\pi$ an irreducible cuspidal $\F$-representation $[\pi]_\ell$ of $\GL nF$. For $\sigma_1,\sigma_2 \in \Go nF_{\Q}$, with determinants of finite order, we say $\sigma_1\equiv \sigma_2\pmod\ell$ if $[\sigma_1]_\ell^\roman{ss} = [\sigma_2]_\ell^\roman{ss}$, and similarly on the other side. We prove: 
\proclaim{Main Theorem} 
If $\sigma_1, \sigma_2 \in \Go nF_{\Q}$ have determinant characters $\det\sigma_1$, $\det\sigma_2$ of finite order, then 
$$ 
\sigma_1\equiv\sigma_2 \pmod \ell \quad \Longleftrightarrow \quad \upr L\sigma_1 \equiv \upr L\sigma_2 \pmod \ell. 
$$ 
\endproclaim 
We use the classification of irreducible cuspidal representations of $\GL nF$ in terms of simple types, as in \cite{8}: although written there in terms of $\Bbb C$-representations, it applies unchanged to $\Q$-representations. 
Thanks to fundamental results of James \cite{17} concerning cuspidal modular representations of finite general linear groups, the relevant classes of simple type behave extraordinarily well under reduction modulo $\ell$. Only the elementary techniques available in \cite{9} or \cite{19}, for example, are necessary to analyze them completely. Another elementary argument relates the congruence behaviour of types to that of cuspidal representations of $\GL nF$. The Comparison Theorem of \cite{6} then implies the desired result for the Langlands correspondence. To reflect this, we have structured the paper to give a brief introduction to the relevant results of \cite{6} alongside an essentially complete and self-contained account of reduction modulo $\ell$. 
\par
Some remarks are in order. First, the map $\sigma\mapsto \upr L\sigma$ on $\Q$-representations does depend on the choice of a field isomorphism $\alpha:\Bbb C\to \Q$. If we replace $\alpha$ by another isomorphism $\beta:\Bbb C\to \Q$, then $\upr L\sigma$ changes to $\upr L(\chi\otimes\sigma)$, where $\chi$ is an unramified character of $\scr W_F$, of order $\le2$, depending only on $F$, $n$ and $\alpha^{-1}\beta$ \cite{13} (7.4). The assertion of the theorem is therefore independent of the choice of $\alpha$, and we  may regard it as an arithmetic property of the complex Langlands correspondence. Next, the techniques and the result hold equally for representations $\sigma\in \Go nF_{\Q}$ such that $\det\sigma$ takes its values in $\Z^\times$, where $\Z$ is the integral closure in $\Q$ of the ring $\Bbb Z_\ell$ of $\ell$-adic integers. However, the reduction from there to the finite case may be achieved by the trivial step of tensoring with an unramified character: we impose the restriction only because it facilitates direct appeal to standard results from finite group theory. 
\par 
The ``modulo $\ell$'' Langlands correspondence already has a substantial literature. Of particular relevance here are the works \cite{20}, \cite{21} of M.-F. Vign\'eras, laying the foundations of a theory of smooth $\F$-representations of $p$-adic reductive groups, and her paper \cite{22} establishing the existence of a semisimple version of the correspondence. More recently, J.-F. Dat has posted a pre-print \cite{10} in which he shows that the standard geometric model, realizing the Langlands correspondence for $\Q$-representations, admits an $\ell$-integral model. On reduction modulo $\ell$, this yields the Langlands correspondence for $\F$-representations. 
\par 
The results of \cite{10} are certainly broader than those here, although the core is the same. However, we feel that the brevity and simplicity of our arguments, presenting the $\F$-correspondence as an easy consequence of an explicit complex result, may be found appealing. While the arguments presented here are short and straightforward,  the proof in its totality cannot be described as elementary or local in nature. It relies via \cite{6} on the standard treatment of the Langlands correspondence \cite{11}, \cite{12}, \cite{18}, along with the classification theory of \cite{8}, and its relations with automorphic induction \cite{14}, \cite{15}, \cite{16}, as developed in \cite{1}, \cite{2} and \cite{5}. 
\head \Smc 
1. Representations of finite groups 
\endhead 
This section is preparatory in nature. We fix a prime number $\ell$ and recall the definition of the classical {\it decomposition map\/} modulo $\ell$. This relates the representations of a finite group over fields of characteristic zero to its representations over fields of characteristic $\ell$. Standard expositions (of which we use \cite{9} and \cite{19}) start from representations over local or global fields, but this will not be convenient for us. We therefore outline the translation to a simpler framework.  We append a couple of specific results used repeatedly later in the paper. 
\subhead 
1.1 
\endsubhead 
Let $\ell$ be a prime number. Let $\Q/\Bbb Q_\ell$ be an algebraic closure of the field $\Bbb Q_\ell$ of $\ell$-adic numbers. Let $\Z$ be the integral closure of $\Bbb Z_\ell$ in $\Q$. Thus $\Z$ is a local ring but, we recall, it is not Noetherian and its unique maximal ideal $\pl$ is not principal. Its residue class field $\Z/\pl$ provides an algebraic closure $\F$ of the field $\Bbb F_\ell$ of $\ell$ elements.  
\proclaim{Definition 1} 
Let $V$ be a finite-dimensional $\Q$-vector space. A \rom{$\Z$-lattice in $V$} is a finitely generated $\Z$-submodule $L$ of $V$ which spans $V$ over $\Q$. 
\endproclaim 
\proclaim{Proposition} 
Let $V$ be a $\Q$-vector space of finite dimension $m$, and let $L$ be a $\Z$-submodule of $V$. The module $L$ is a $\Z$-lattice in $V$ if and only if it is generated over $\Z$ by a $\Q$-basis of $V$. In particular, any $\Z$-lattice in $V$ is $\Z$-free of rank $m$. 
\endproclaim 
\demo{Proof} 
One implication is trivial. To prove the converse, let $L$ be a $\Z$-lattice in $V$. Since $L$ is finitely generated over $\Z$, we may choose a $\Z$-generating set $\{x_1,x_2,\dots,x_r\}$ of $L$, of minimal cardinality. This set must surely span $V$ over $\Q$. Suppose, for a contradiction, that it is linearly dependent over $\Q$, say $\sum_{i=1}^ra_ix_i = 0$, for elements $a_i\in \Q$, not all zero. There is a finite field extension $E/\Bbb Q_\ell$ containing all $a_i$. Scaling by an element of $E^\times$, we may assume that all $a_i$ lie in the discrete valuation ring $\frak o_E$ in $E$, and that at least one of them, say $a_1$, is a unit in $\frak o_E$.  Thus $x_1$ lies in the $\Z$-module generated by $\{x_2,x_3,\dots,x_r\}$, contrary to hypothesis. It follows that $\{x_1,x_2,\dots,x_r\}$ is a $\Q$-basis of $V$ (and a $\Z$-basis of $L$).  \qed 
\enddemo 
\proclaim{Corollary} 
Let $V$ be a finite-dimensional $\Q$-vector space, let $U$ be a subspace of $V$ and let $L$ be a $\Z$-lattice in $V$. The group $L\cap U$ is then a $\Z$-lattice in $U$. 
\endproclaim 
\demo{Proof} 
The image $L/L\cap U$ of $L$ in $V/U$ is a finitely generated $\Z$-module, spanning $V/U$ over $\Z$. It is, by the proposition, $\Z$-free of finite rank. Thus $L = L\cap U\oplus L'$, for a finitely generated free $\Z$-submodule $L'$ of $L$, isomorphic to $L/L\cap U$. In particular, $L\cap U \cong L/L'$, so $L\cap U$ is finitely generated over $\Z$. If $W$ is the $\Q$-span of $L'$, then $V = U\oplus W$ and so $L\cap U$ spans $U$. Thus $L\cap U$ is a $\Z$-lattice in $U$, as required. \qed 
\enddemo 
\proclaim{Definition 2} 
Let $V$ be a $\Q$-vector space, of possibly infinite dimension. A $\Z$-submodule $L$ of $V$ is a $\Z$-lattice in $V$ if $L\cap U$ is a $\Z$-lattice in $U$, for every finite-dimensional subspace $U$ of $V$. 
\endproclaim 
The corollary implies that the two definitions are consistent. 
\subhead 
1.2 
\endsubhead 
Let $G$ be a {\it finite\/} group. If $k$ is a field, we denote by $\frak G_0(kG)$ the Grothendieck group of finite-dimensional representations of $G$ over $k$, formed relative to exact sequences. Let $\frak R$ be the set of isomorphism classes of irreducible $k$-representations of $G$. The group $\frak G_0(kG)$ is then free abelian on $\frak R$. We may identify the set of formal sums $\sum_{R\in \frak R} n_RR$, in which all $n_R$ are non-negative, with the set of isomorphism classes of finite-dimensional {\it semisimple\/} $k$-representations of $G$. 
\par 
We consider the category of finite-dimensional representations of $G$ over the field $\Q$. The fields $\Q$, $\Bbb C$ are isomorphic, and the choice of an isomorphism $\Bbb C \to \Q$ induces an equivalence between the categories of finite-dimensional representations of $G$ over the two fields. We may so translate the basic theory of complex representations to the new context unchanged. 
\par 
Let $(\rho,U)$ be a representation of $G$ over $\Q$, of finite dimension $m$. A {\it $\Z G$-lattice in\/} $U$ is a $\Z$-lattice in $U$ which is stable under the action of $\rho(G)$. If $L_0$ is a $\Z$-lattice in $U$, then the $\Z$-module $\sum_{g\in G} \rho(g)L_0$ is $G$-stable and it is a $\Z$-lattice. Thus $(\rho,U)$ admits a $\Z G$-lattice. 
\par
Let $L$ be a $\Z G$-lattice in $U$. By 1.1 Proposition, the quotient space $\wt L = L/\pl L$ is a vector space over $\F$, of dimension $m$. The representation $\rho$ induces a representation $\rho_L$ of $G$ on $\wt L$. Let $\rho_L^\roman{ss}$ denote the image of $\rho_L$ in the Grothendieck group $\frak G_0(\F G)$. 
\proclaim{Proposition} 
The class $\rho^\roman{ss}_L$ depends only on the isomorphism class of $(\rho,U)$, and not on the choice of\/ $\Z G$-lattice $L$ in $U$. In particular, if $\rho_L$ is irreducible, then its isomorphism class is determined by that of $(\rho,U)$, independently of the choice of $L$.  
\endproclaim 
\demo{Proof} 
We perform a preliminary reduction. Let $K/\Bbb Q_\ell$ be a finite field extension, $K\subset \Q$, such that $K$ contains a primitive $|G|$-th root of unity. Extension of scalars then induces an isomorphism $\frak G_0(KG) \to \frak G_0(\Q G)$ (\cite{19} Th\'eor\`eme 24, Corollaire). In particular, there is an irreducible $K$-representation $(\rho_0,U_0)$ of $G$ such that $(\rho,U)$ is isomorphic to the representation of $G$ on $\Q\otimes_KU_0$ induced by $\rho_0$. We henceforward identify $U$ with $\Q \otimes_K U_0$. 
\par 
Let $L$ be a $\Z G$-lattice in $U$. Let $\frak o_K$ be the discrete valuation ring in $K$ and let $L_0$ be a $G$-stable $\frak o_K$-lattice in $U_0$. We set $\overline L_0 = \Z\otimes_{\frak o_K}L_0$. Thus $\overline L_0$ is a $G$-stable $\Z$-lattice in $U$. Expressing a $\Z$-basis of $L$ in terms of an $\frak o_K$-basis of $L_0$, we see: 
\proclaim{Lemma} 
Let $L$ be a $\Z G$-lattice in $U$. There is a finite field extension $K'/K$ and an $\frak o_{K'}G$-lattice $L'_0$ in $K'\otimes_K U_0$ such that $L = \Z \otimes_{\frak o_{K'}} L'_0$. 
\endproclaim 
Let $L_1$, $L_2$ be $\Z G$-lattices in $U$. According to the lemma, we may choose the field $K$ and $\frak o_KG$-lattices $L_{i0}$ in $U_0$ such that $L_i = \Z\otimes_{\frak o_K} L_{i0}$, $i=1,2$. Let $\frak p_K = \pl\cap K$ be the maximal ideal of $\frak o_K$ and write $\Bbbk_K = \frak o_K/\frak p_K$. By \cite{9} (16.16) or \cite{19} Th\'eor\`eme 32, the $\Bbbk_K$-representations $\rho_{L_{i0}}$ of $G$ on $L_{i0}/\frak p_KL_{i0}$, induced by $\rho_0$, define the same element of $\frak G_0(\Bbbk_KG)$. The representation $\rho_{L_i}$ is obtained from $\rho_{L_{i0}}$ by extension of scalars from $\Bbbk_K$ to $\F$. Therefore $\rho_{L_1}^\roman{ss}= \rho_{L_2}^\roman{ss}$, as required. The final assertion follows. \qed
\enddemo 
The proposition shows that $\rho\mapsto \rho_L^\roman{ss}$ induces a well-defined map $\frak G_0(\Q G) \to \frak G_0(\F G)$. We denote it 
$$ 
\aligned 
\frak G_0(\Q G) &\longrightarrow \frak G_0(\F G), \\ 
\rho &\longmapsto [\rho]_\ell^\roman{ss}. 
\endaligned 
\tag 1.2.1 
$$
\remark{Remark} 
The map (1.2.1) is the {\it decomposition map\/} of the literature, in which it is usually denoted $d_\ell$. We often think of $[\rho]_\ell^\roman{ss}$ as the isomorphism class of the semisimplification of the representation $\rho_L$, for some lattice $L$. 
\endremark 
\subhead 
1.3 
\endsubhead 
We continue this general discussion with a simple remark which will be useful later. 
\proclaim{Proposition} 
Let $(\rho,V)$ be an irreducible $\Q$-representation of the finite group $G$, and suppose that $[\rho]_\ell^\roman{ss}$ is irreducible. If $L$, $L'$ are $\Z G$-lattices in $V$, there exists $a\in \Q^\times$ such that $L' = aL$. 
\endproclaim 
\demo{Proof} 
As in the proof of 1.1 Proposition, there exists $a\in \Q^\times$ such that $L\supset aL' \not\subset \frak P_\ell L$. Replacing $L'$ by $aL'$, we may assume $L\supset L'\not\subset \frak P_\ell L$. The image of $L'$ in $L/\frak P_\ell L$ is then a non-zero $G$-subspace of the irreducible $\F G$-space $L/\frak P_\ell L$. We deduce that $L = L'{+}\frak P_\ell L$. 
\par 
We show that this condition implies $L' = L$, whence the result will follow. We choose a $\Z$-basis $\scr B$ of $L$ and a $\Z$-basis $\scr B'$ of $L'$. The ``transition matrix'', relating $\scr B'$ to $\scr B$, then has coefficients in $\Z$, and its reduction modulo $\frak P_\ell$ is invertible over $\F$. The matrix is therefore invertible over $\Z$, whence $L' = L$ as required. \qed 
\enddemo 
\subhead 
1.4 
\endsubhead 
We record some special properties of $p$-groups, where $p$ is a prime number other than $\ell$. 
\proclaim{Proposition} 
Let $p$ be a prime number, $p\neq \ell$. Let $G$ be a finite $p$-group. Any $\F$-representation of $G$ is semisimple. 
\par 
If $\rho_1$, $\rho_2$ are irreducible $\Q$-representations of $G$, then  
\roster 
\item 
the $\F$-representations $[\rho_1]_\ell^\roman{ss}$, $[\rho_2]_\ell^\roman{ss}$ are irreducible, and  
\item 
$[\rho_1]_\ell^\roman{ss} = [\rho_2]_\ell^\roman{ss}$ if and only if $\rho_1\cong \rho_2$. 
\endroster 
\endproclaim 
\demo{Proof} The first assertion is Maschke's Theorem. The others are equivalent to \cite{19} Proposition 43. \qed 
\enddemo 
\subhead 
1.5 
\endsubhead 
For our next result, we are given a finite group $G$ and a subgroup $H$. We define the induction functors $\Ind_H^G$, $\cind_H^G$ following exactly the formulation for complex representations of locally profinite groups. In this situation, they are the same so we get two versions of Frobenius Reciprocity: if $\rho$ is an $\F$-representation of $G$, and $\tau$ an $\F$-representation of $H$, then 
$$
\align 
\Hom G\rho{\Ind_H^G\,\tau} &\cong \Hom H \rho\tau, \\ 
\Hom G{\Ind_H^G\,\tau}\rho &\cong \Hom H\tau\rho. 
\endalign 
$$ 
\proclaim{Proposition} 
Let $p$ be a prime number, $p\neq \ell$. Let $G$ be a finite group, and let $P$ be a normal $p$-subgroup of $G$. Let $\alpha$ be an irreducible $\F$-representation of $P$, and let $H$ be the group of $g\in G$ such that $\alpha^g \cong\alpha$. Suppose there exists a representation $\rho$ of\/ $H$ such that $\rho|_{P} \cong \alpha$. 
\roster 
\item 
Let $\tau$ be an irreducible $\F$-representation of $H$, trivial on $P$. The representation 
$$
\nu_\tau = \Ind^G_H\,\rho\otimes \tau 
$$
is irreducible. 
\item 
The map $\tau\mapsto \nu_\tau$ is a bijection between the set of isomorphism classes of irreducible $\F$-representations $\tau$ of\/ $H$ such that $\tau|_P$ is trivial, and the set of isomorphism classes of irreducible representations $\nu$ of $G$ such that $\Hom P\alpha\nu \neq 0$. 
\endroster 
\endproclaim 
\demo{Proof} 
In (1), set $\dim\tau = d_\tau$, and consider the semisimple representation $\nu_\tau|_P$. Since $P$ is normal in $G$, it is the direct sum of representations $\alpha^g$, $g\in H\backslash G$, each occurring with multiplicity $d_\tau$. If $\check\rho$ denotes the contragredient of $\rho$, it follows that the space of $P$-fixed points in $\check\rho\otimes \nu_\tau$ has dimension $d_\tau$. This space, moreover, carries a representation $\tau'$ of $H$, trivial on $P$. To identify the representation $\tau'$, we observe that 
$$ 
\Hom H\tau{\tau'} = 
\Hom H\tau{\check\rho\otimes \nu_\tau} \cong \Hom H{\rho\otimes\tau}{\nu_\tau} \cong \Hom G {\nu_\tau}{\nu_\tau} \neq 0. 
$$
Since $\dim\tau' = d_\tau = \dim\tau$, we conclude that $\tau'$ is equivalent to $\tau$. We further deduce that $\Hom G {\nu_\tau} {\nu_\tau} = \F$. 
\par 
Next, let $\mu$ be an irreducible $G$ sub-representation of $\nu_\tau$. The restriction of $\mu$ to $P$ is again a direct sum of representations $\alpha^g$, $g\in H\backslash G$, each occurring with the same multiplicity. The space of $P$-fixed points in the representation $\check\rho\otimes\mu$ is therefore nonzero, and carries a representation $\mu_0$ of $H$, trivial on $P$. This representation $\mu_0$ appears naturally as a sub-representation of $\tau'$ and, since $\tau ' \cong \tau$ is irreducible, we get $\mu_0 \cong \tau$. Therefore 
$$
\Hom G{\nu_\tau}\mu \cong \Hom H \tau{\check\rho\otimes \mu} \cong \F. 
$$ 
The space $\Hom G{\nu_\tau}\mu$ is a subspace of $\Hom G{\nu_\tau}{\nu_\tau}$, while 
$$
\Hom G{\nu_\tau}{\nu_\tau} \cong \Hom H \tau{\check\rho\otimes \nu_\tau} \cong \F. 
$$ 
The identity endomorphism of $\nu_\tau$ therefore lies in $\Hom G{\nu_\tau}\mu$, whence $\mu = \nu_\tau$. In particular, $\nu_\tau$ is irreducible, as required for (1). Since $\tau$ appears as the natural representation of $H$ on the space of $P$-fixed points in $\check\rho\otimes \nu_\tau$, we have also proved (2). \qed 
\enddemo 
\remark{Remark} 
In the proposition, the hypothesis $P\triangleleft G$ is essential. Given that, the result and its proof are valid for representations of $G$ over any algebraically closed field of characteristic other than $p$. 
\endremark 
\head\Smc 
2. Linear groups over local fields 
\endhead 
Again, $\ell$ is a prime number. From now on, we work with a non-Archimedean local field $F$ with residue class field $\Bbbk_F$ of characteristic $p\neq \ell$. We mildly generalize the machinery of \S1 to admissible $\Q$-representations of the locally profinite group $G = \GL nF$. 
\subhead 
2.1 
\endsubhead 
Let $(\pi,V)$ be an {\it admissible\/} $\Q$-representation of $G$, and let $L$ be a $\Z$-lattice in $V$, in the sense of 1.1 Definition 2. In particular, if $K$ is a compact open subgroup of $G$, the space $V^K$ of $\pi(K)$-fixed points in $V$ has finite dimension and $L\cap V^K$ is a $\Z$-lattice in $V^K$. 
\par 
We say that $L$ is a $\Z G$-lattice if $\pi(g)L = L$, for all $g\in G$. In such a case, the lattice $L\cap V^K$ is the set $L^K$ of $\pi(K)$-fixed points in $L$. More generally, let $\rho$ be an irreducible, smooth $\Q$-representation of $K$. The isotypic subspace $V^\rho$ of $V$ has finite dimension, and so $L^\rho = L\cap V^\rho$ is a $\Z K$-lattice in $V^\rho$. 
\proclaim{Proposition} 
Let $(\pi,V)$ be an admissible $\Q$-representation of $G$, and let $L$ be a $\Z G$-lattice in $V$. Let $K$ be an open pro-$p$ subgroup of $G$, and let $\widehat K$ denote the set of equivalence classes of irreducible smooth $\Q$-representations of $K$. For $\rho\in \widehat K$, the set $L^\rho$ is a $\Z K$-lattice in $V^\rho$ and $L = \bigoplus_{\rho\in \widehat K} L^\rho$. 
\endproclaim 
\demo{Proof} 
The first assertion has already been remarked. Let $K'$ be an open normal subgroup of $K$. We view $V^{K'}$ as providing a representation of the finite group $K/K'$, and then $L^{K'}$ is a $\Z[K/K']$-lattice in $V^{K'}$. Moreover, $V^{K'} = \bigoplus_\rho V^\rho$, where $\rho$ ranges over those elements of $\widehat K$ with $K'\subset \roman{Ker}\,\rho$. For such $\rho$, there is a primitive central idempotent $e_\rho$ of the group algebra $\Q[K/K']$ such that $V^\rho = e_\rho V^{K'}$. The integer $(K{:}K')$ is a power of $p$, so it lies in $\Z^\times$. Therefore $e_\rho \in \Z [K/K']$, whence $L^{K'} = \bigoplus_\rho L^\rho$, with $\rho$ ranging as before. Allowing $K'$ to range over a descending sequence of open normal subgroups of $K$ with trivial intersection, the result follows. \qed 
\enddemo 
\remark{Remark} 
We have not asserted the existence of a $\Z G$-lattice in an arbitrary admissible representation of $G$. In the one case we need, we will give a direct construction. For a more general discussion of lattices, see \cite{20} I, \S9. 
\endremark 
\subhead 
2.2 
\endsubhead 
We need one more simple property. 
\proclaim{Proposition} 
Let $(\pi,V)$ be an admissible $\Q$-representation of $G$, admitting a $\Z G$-lattice $L$. Let $\wt L = L/\frak P_\ell L$ and let $\pi_L$ denote the natural $\F$-representation of $G$ on $\wt L$. Let $K$ be an open pro-$p$ subgroup of $G$. 
\roster 
\item 
If $q:L\to \wt L$ is the quotient map, then $q(L^K) = \wt L^K$. In particular, the representation $\pi_L$ is admissible. 
\item 
Let $\rho$ be an irreducible smooth $\Q$-representation of $K$. If $\tilde\rho$ denotes the irreducible $\F$-representation $[\rho]_\ell^\roman{ss}$, then $q(L^\rho) = \widetilde L^{\tilde\rho}$. 
\endroster 
\endproclaim 
\demo{Proof} 
Let $\tilde v\in L^K$ and choose $v \in L$ such that $q(v) = \tilde v$. There then exists an open normal subgroup $K'$ of $K$ such that $v\in L^{K'}$. The element $v' = (K{:}K')^{-1}\sum_{x\in K/K'} \pi(x)v$ lies in $L^K$ and satisfies $q(v') = \tilde v$. This proves part (1), and part (2) now follows from 2.1 Proposition. \qed 
\enddemo 
\head\Smc 
3. Induced cuspidal representations 
\endhead 
We consider irreducible {\it cuspidal\/} representations of the group $G = \GL nF$ over the field $\Q$. The classification and structure theory for complex cuspidal representations of $G$, laid out in \cite{8} and further developed in the first nine sections of \cite{1}, are purely algebraic in nature: we may apply them unchanged to $\Q$-representations via a field isomorphism $\Q\to \Bbb C$. We investigate the behaviour of these structures under reduction modulo $\ell$. 
\par 
We start by reviewing the background, relying as much as possible on the summary given in \cite{6} (especially sections 2.1 and 2.2). 
\subhead 
3.1 
\endsubhead 
Let $\theta$ be an {\it m-simple character\/} in $G$. Thus, by definition, $\theta$ is either the trivial character of $U^1_\frak m$, where $\frak m$ is a maximal order in $A = \M nF$, or else $\theta$ is a simple character attached to a simple stratum $[\frak a,\beta]$ in $A$, in which $\frak a$ is maximal among the hereditary $\frak o_F$-orders in $A$ that are stable under conjugation by the group $F[\beta]^\times$ of non-zero elements of the field $F[\beta]$. In particular, $\theta$ is a character of the compact open subgroup $H^1_\theta = H^1(\beta,\frak a)$ of $G$, and $H^1_\theta$ is a pro-$p$ group. 
\par 
In all cases, we let $\bk J_\theta$ be the $G$-normalizer of $\theta$, we let $J^0_\theta$ be the unique maximal compact subgroup of $\bk J_\theta$ and $J^1_\theta$ the pro-$p$ radical of $J^0_\theta$. An {\it extended maximal simple type\/} over $\theta$ is an irreducible representation of $\bk J_\theta$ which contains $\theta$ and is intertwined only by elements of $\bk J_\theta$. Let $\scr T(\theta)$ be the set of equivalence classes of extended maximal simple types over $\theta$. We remark ({\it cf\.} 3.2 below) that a representation $\vL$ of $\bk J_\theta$ lies in $\scr T(\theta)$ if and only if $(J^0_\theta, \vL|_{J^0_\theta})$ is a {\it maximal simple type in $G$}, in the sense of \cite{8}. We summarize the main points. 
\proclaim{Proposition} 
\roster 
\item Let $(\pi,V)$ be an irreducible cuspidal representation of $G$. 
\itemitem{\rm (a)} The representation $\pi$ contains a simple character $\theta$. The character $\theta$ is m-simple and is unique up to $G$-conjugacy. 
\itemitem{\rm (b)} The natural representation $\vL_\pi$ of $\bk J_\theta$ on $V^\theta$ is irreducible and lies in $\scr T(\theta)$. In particular, $\pi \cong \cind_{\bk J_\theta}^G\,\vL_\pi$. 
\item Let $\theta$ be an m-simple character in $G$. The map $\vL \mapsto \cind_{\bk J_\theta}^G\,\vL$, $\vL\in \scr T(\theta)$ is a bijection between $\scr T(\theta)$ and the set of equivalence classes of irreducible cuspidal representations of $G$ which contain $\theta$. 
\endroster 
\endproclaim 
\subhead 
3.2 
\endsubhead 
Let $\theta$ be an m-simple character in $G$. We describe the elements of $\scr T(\theta)$ more explicitly. 
\par 
Suppose first that $\theta$ is the trivial character of $U^1_\frak m$, for a maximal $\frak o_F$-order $\frak m$ in $A$. In this case, $\bk J_\theta = F^\times J^0_\theta$, $J^0_\theta = U_\frak m$ and $J^1_\theta = H^1_\theta = U^1_\frak m$. We have $J^0_\theta/J^1_\theta = U_\frak m/U^1_\frak m \cong \GL n{\Bbbk_F}$. A representation $\vL$ of $\bk J_\theta$ then lies in $\scr T(\theta)$ if and only if $\vL|_{J^0_\theta}$ is the inflation of an irreducible {\it cuspidal\/} representation of $J^0_\theta/J^1_\theta \cong \GL n{\Bbbk_F}$. Any two maximal orders in $A$ are $G$-conjugate, so the set of $G$-conjugacy classes of elements of $\scr T(\theta)$ effectively depends only on the dimension $n$. 
\par 
We therefore assume $\theta$ is non-trivial, attached to a simple stratum $[\frak a,\beta]$. Denoting the field $F[\beta]$ by $P$, we here have $\bk J_\theta = P^\times J^0_\theta$. If $C$ denotes the $A$-centralizer of $P$, the set of elements of $G$ which intertwine $\theta$ is $I_G(\theta) = J^1_\theta C^\times J^1_\theta$. 
\par 
We define two sets of irreducible representations of $\bk J_\theta$. Let $\eta_\theta$ be the unique irreducible representation of $J^1_\theta$ containing $\theta$. Let $\scr H(\theta)$ denote the set of equivalence classes of representations $\kappa$ of $\bk J_\theta$ such that $\kappa|_{J^1_\theta} \cong \eta_\theta$ and such that $\kappa$ is intertwined by every element of $I_G(\theta)$. 
\par 
Let $\frak c = \frak a\cap C$. Since $\theta$ is m-simple, $\frak c$ is a maximal $\frak o_P$-order in $C$ and so isomorphic to $\M m{\frak o_P}$, $n=m[P{:}F]$. We have $J^0_\theta = U_\frak cJ^1_\theta$, and $U_\frak c\cap J^1_\theta = U^1_\frak c$. We define $\scr T_0(\theta)$ to be the set of equivalence classes of irreducible representations $\lambda$ of $\bk J_\theta$ such that $\lambda$ is trivial on $J^1_\theta$ and $\lambda|_{J^0_\theta}$ is inflated from an irreducible cuspidal representation of $J^0_\theta/J^1_\theta \cong U_\frak c/U^1_\frak c$. 
\remark{Remark} 
If we write $\theta_\frak c$ for the trivial character of $U^1_\frak c$, then $\theta_\frak c$ is a trivial m-simple character in $C^\times$, and the map $\lambda\mapsto \lambda|_{P^\times U_\frak c}$ gives a bijection $\scr T_0(\theta) \to \scr T(\theta_\frak c)$. 
\endremark 
\proclaim{Proposition} 
Let $\kappa\in \scr H(\theta)$, $\lambda \in \scr T_0(\theta)$. The representation $\kappa\otimes \lambda$ lies in $\scr T(\theta)$, and the map $\scr T_0(\theta) \to \scr T(\theta)$, $\lambda \mapsto \kappa\otimes \lambda$, is a bijection. 
\endproclaim 
\demo{Proof} See 3.6 Proposition of \cite{6}. \qed 
\enddemo 
\subhead 
3.3 
\endsubhead 
We use these facts to investigate the reduction properties of cuspidal representations. 
\par 
Let $\theta$ be an m-simple character in $G$, and let $\vL\in \scr T(\theta)$. The restriction of $\vL$ to $F^\times$ is a multiple of a character $\omega_\vL$. The group $\vL(\bk J_\theta)$ is then finite if and only if $\omega_\vL(F^\times)$ is finite. So, when $\omega_\vL$ has finite order, the representation $\vL$ admits a $\Z \bk J_\theta$-lattice. 
\proclaim{Proposition} 
Let $(\bk J_\theta,\vL,W)\in \scr T(\theta)$, suppose that $\omega_\vL$ has finite order, and let $L$ be a $\Z\bk J_\theta$-lattice in $W$. 
\roster 
\item 
The $\F$-representation $\vL_L$ of $\bk J_\theta$ on $\wt L = L/\pl L$ is irreducible. 
\item 
The isomorphism class of $\vL_L$ is independent of the choice of $L$. 
\item 
If $L'$ is a $\Z \bk J_\theta$-lattice in $W$, there exists $a\in \Q^\times$ such that $L' = aL$. 
\endroster 
\endproclaim 
\demo{Proof} 
We observe that (1) implies (2) (by 1.2 Proposition) and (2) implies (3) (by 1.3 Proposition). 
\par 
To prove (1), we consider first the case where the m-simple character $\theta$ is trivial. Thus, we may assume, $\theta$ is the trivial character of $U^1_\frak m$, where $\frak m = \M n{\frak o_F}$. The group $\bk J_\theta$ is $F^\times \GL n{\frak o_F}$ and, by definition, the restriction of $\vL$ to $\GL n{\frak o_F}$ is the inflation of an irreducible cuspidal representation of $\GL n{\Bbbk_F}$. The assertion (1) is thus an instance of \cite{17} Theorem 3.6. 
\par 
We therefore assume $\theta$ to be non-trivial, attached to a simple stratum $[\frak a,\beta]$, as in 3.2. In particular, $\vL = \kappa\otimes \lambda$, for  representations $\kappa\in \scr H(\theta)$ and $\lambda\in \scr T_0(\theta)$. We may choose $\kappa$ such that $\kappa(F^\times)$ is finite \cite{6} 3.2 Corollary; in this case, $\lambda(F^\times)$ is also finite.  
\par  
Let $\lambda$ (resp\. $\kappa$) act on the vector space $W_1$ (resp\. $W_2$), and let $L_i$ be a $\Z\bk J_\theta$-lattice in $W_i$. As in the first case of the proof, the representation $\lambda_{L_1}$ is irreducible. Let $\eta_\theta$ be the unique irreducible representation of $J^1_\theta$ containing $\theta$. Since $J^1_\theta$ is a pro-$p$ group, the representation $(\eta_\theta)_{L_2}$ is irreducible (1.4 Proposition). As $\kappa|_{J^1_\theta} \cong \eta_\theta$, the representation $\kappa_{L_2}$ is also irreducible. The representation $\vL$ acts on $W = W_1\otimes W_2$, and $L = L_1\otimes L_2$ provides a $\Z\bk J_\theta$-lattice in $W$. The representation $\vL_L = \lambda_{L_1}\otimes \kappa_{L_2}$ is irreducible, by 1.5 Proposition. \qed 
\enddemo 
\subhead 
3.4 
\endsubhead 
We continue in the same situation, with an m-simple character $\theta$ in $G = \GL nF$, and an extended maximal simple type $\vL\in \scr T(\theta)$ such that $\omega_\vL$ has finite order. The representation $(\pi,V) = \cind_{\bk J_\theta}^G\,\vL$ is irreducible, cuspidal and the central character $\omega_\pi = \omega_\vL$ has finite order. We let $\vL$ act on the $\Q$-space $W$, and choose a $\Z\bk J_\theta$-lattice $L$ in $W$. The space $V$ consists of all functions $f:G\to W$ which are compactly supported modulo the centre of $G$ and satisfy $f(jg) = \vL(j)f(g)$, for $j\in \bk J_\theta$, $g\in G$. The group $G$ acts on $V$ by right translation. We define 
$$
L_* = \{f\in V: f(G)\subset L\}. 
\tag 3.4.1
$$ 
Surely, $L_*$ is a $\pi(G)$-stable $\Z$-submodule of $V$. 
\proclaim{Lemma} 
The $\Z$-module $L_*$ is a $\Z G$-lattice in $V$. 
\endproclaim 
\demo{Proof} 
It is enough to show that $L_*^K$ is a $\Z$-lattice in $V^K$, for any open, pro-$p$ subgroup $K$ of $G$. We follow the procedure of \cite{4} 3.5 to write down a basis of the (finite-dimensional) space $V^K$. First, we choose a set $X$ of coset representatives for $\bk J_\theta\backslash G/K$. For each $x\in X$, we choose a basis $\scr B_x$ of the space of $xKx^{-1}\cap \bk J_\theta$-fixed points in $W$. For each $x\in X$ and each $w\in \scr B_x$, there is a unique function $f_{x,w}\in V^K$ with support $\bk J_\theta xK$ and such that $f_{x,w}(x) = w$. The set 
$$
\scr B(K) = \{f_{x,w}:x\in X, w\in \scr B_x\} 
$$
is then a basis of $V^K$. For each $x\in X$, we may choose $\scr B_x$ to be a $\Z$-basis of the lattice $L^{xKx^{-1}\cap \bk J_\theta}$. The corresponding set $\scr B(K)$ is then a $\Z$-basis of $L_*^K$. \qed 
\enddemo 
We form the $\F$-representation $(\pi_{L_*}, L_*/\pl L_*)$ of $G$: this is admissible by 2.2 Proposition. 
\proclaim{Proposition} 
The representation $(\pi_{L_*}, L_*/\pl L_*)$ is irreducible, and equivalent to the representation $\cind_{\bk J_\theta}^G\, \vL_L$. 
\endproclaim 
\demo{Proof} 
Let $f\in L_*$; we define a function $f_L:G\to \wt L$ by composing $f$ with the canonical map $L\to \wt L$. If $f\in \pl L_*$, then $f(G)\subset \pl L$ and $f_L = 0$. We thus obtain a $G$-homomorphism $L_*/\pl L_* \to \cind\, \vL_L$. If $K$ is an open pro-$p$ subgroup of $G$, we follow the procedure of the lemma to construct a $\Z$-basis $\scr B(K)$ of $L_*^K$. The set $\{b_L:b\in \scr B(K)\}$ then provides an $\F$-basis of the space $\big(\cind\,\vL_L\big)^K$. It follows that the map 
$$
\align 
L_* &\longrightarrow \cind_{\bk J_\theta}^G\,\vL_L, \\ f &\longmapsto f_L, \endalign 
$$ 
induces an $\F G$-isomorphism $L_*/\pl L_* \cong \cind \,\vL_L$. 
\par
We identify $L$ with the sublattice of $L_*$ consisting of functions with support $\bk J_\theta$. The lattice $L$ is also the $\theta$-isotypic sublattice $L_*^\theta = L_*\cap V^\theta$ of $L_*$. The character $\theta$ has finite order, which is a power of $p$. Let $\mu_p(\Q)$ denote the group of roots of unity in $\Q$ of $p$-power order. Thus $\mu_p(\Q) \subset \Z^\times$ and, under the reduction map $\Z^\times \to \F^\times$, $\mu_p(\Q)$ is mapped isomorphically to the group $\mu_p(\F)$ of $p$-power roots of unity in $\F^\times$. Consequently, there is no need to distinguish between $\theta$ and $[\theta]_\ell^\roman{ss}$,  its reduction mod\. $\ell$. With this convention, 2.2 Proposition implies that the canonical map $L_* \to L_*/\pl L_*$ identifies $L_*^\theta/\pl L_*^\theta$ with $(L_*/\pl L_*)^\theta$, that is, with $\wt L \subset L_*/\pl L_*$. 
\par 
Let $U$ be a non-zero $G$-subspace of $\cind\,\vL_L$. We apply Frobenius Reciprocity, 
$$
0\neq \Hom GU{\cind\,\vL_L} \subset \Hom GU{\Ind\,\vL_L} \cong \Hom{\bk J_\theta} U{\vL_L}. 
$$ 
In particular, $U^\theta\neq 0$. However, the isotypic space $(\cind\,\vL_L)^\theta = \wt L$ is irreducible over $\bk J_\theta$, so $U\supset \wt L$. Since $\wt L$ generates $\cind\,\vL_L$ over $G$, we have $U = \cind\,\vL_L$, as required. \qed 
\enddemo 
\remark{Remark} 
The standard argument, as in for example \cite{4} 11.4, shows that the irreducible $\F$-representation $\pi_{L_*}$ is {\it cuspidal,} in that its matrix coefficients are compactly supported modulo $F^\times$. As in \cite{20} II 2.3, 2.7, this is equivalent to $\pi_{L_*}$ not being equivalent to a sub-representation of any proper parabolically induced representation. 
\endremark
\subhead 
3.5 
\endsubhead 
We consider general $\Z G$-lattices in an irreducible cuspidal representation. 
\proclaim{Proposition} 
Let $(\pi,V)$ be an irreducible, cuspidal $\Q$-representation of $G$, such that $\omega_\pi$ has finite order. 
\roster 
\item 
The representation $(\pi,V)$ admits a $\Z G$-lattice. 
\item 
Let $M$ be a $\Z G$-lattice in $(\pi,V)$. The natural representation $\pi_M$ of $G$ on $M/\pl M$ is irreducible. The isomorphism class of $\pi_M$ depends only on that of $\pi$, and not on $M$. 
\endroster 
\endproclaim 
\demo{Proof} 
The representation $(\pi,V)$ contains an m-simple character $\theta$. Let $\vL$ denote the natural representation of $\bk J_\theta$ on $V^\theta$, so that $\pi \cong \cind_{\bk J_\theta}^G\,\vL$. As in 3.4 Proposition, we choose a $\Z \bk J_\theta$-lattice $L$ in $V^\theta$. According to that result, the set $L_*$ of (3.4.1) provides a $\Z G$-lattice in $V$, as required for part (1). 
\par 
For part (2), let $M$ be a $\Z G$-lattice in $V$. The $\Z$-module $M^\theta = M\cap V^\theta$ provides a $\Z \bk J_\theta$-lattice in $\vL$ (2.1 Proposition). It follows from 1.3 Proposition that $M^\theta = aL$, for some $a\in \Q^\times$. Scaling, therefore, we may assume that $M^\theta = L = L_*^\theta$. Since $L_*$ is generated by $L$ as $\Z G$-lattice, we have $L_* \subset M$. 
\par 
Suppose, for a contradiction, that $M\neq L_*$. There then exists an open, pro-$p$ subgroup $K$ of $\bk J_\theta$ such that $M^K \varsupsetneq L_*^K$. We choose $x\in \Q^\times$ such that, if $M_1$ denotes the lattice $xM$, then $M^K_1$ is contained in $L_*^K$ but not in $\pl L_*^K$. Let $M_2$ be the $\Z G$-lattice generated by $M^K_1$. Thus $\pl L_* \not\supset M_2\subset L_*$, but $M_2^\theta = M_1^\theta$ is contained in $\pl L$. The image of $M_2$ in $L_*/\pl L_*$ is then a proper $G$-subspace, which is impossible. This proves (2). \qed 
\enddemo 
\remark{Remark} 
We have also shown that the cuspidal representation $(\pi,V)$ admits only one $\Z G$-lattice, up to scaling by elements of $\Q^\times$. 
\endremark 
\subhead 
3.6 
\endsubhead 
Let $\Ao nF$ denote the set of equivalence classes of irreducible, cuspidal, $\Q$-representations of $G = \GL nF$, and let $\Ao nF^{\text{fin}}$ be the set of $(\pi,V)\in \Ao nF$ such that $\omega_\pi$ has finite order. For $(\pi,V)\in \Ao nF^{\text{fin}}$, we define $[\pi]_\ell$ to be the isomorphism class of the natural representation $\pi_M$ of $G$ on $M/\pl M$, for a $\Z G$-lattice $M$ in $V$. As in 3.5 Proposition, $[\pi]_\ell$ is irreducible and independent of the choice of $M$. We may therefore make the following definition.  
\proclaim{Definition} 
Let $\pi_1,\pi_2\in \Ao nF^{\text{fin}}$; say $\pi_1\equiv \pi_2 \pmod\ell$ if $[\pi_1]_\ell = [\pi_2]_\ell$. 
\endproclaim 
\remark{Remark} 
The group $\roman{Aut}\,\Q$ of field automorphisms of $\Q$ acts on $\Ao nF^\roman{fin}$: we denote this action by $\alpha:\pi\mapsto \pi^\alpha$. This action is compatible with $\Z$-structures, and we conclude that, if $\alpha\in \roman{Aut}\,\Q$ and $\pi_i\in \Ao nF^\roman{fin}$, then $\pi_1^\alpha \equiv \pi_2^\alpha \pmod \ell$ if and only if $\pi_1 \equiv \pi_2 \pmod \ell$. Consequently, we can transfer the congruence relation to complex representations as follows. Let $\pi_1$, $\pi_2$ be irreducible, cuspidal $\Bbb C$-representations of $G$, each with central character of finite order. By composing with a field isomorphism $\gamma:\Bbb C \to \Q$, we obtain representations $\pi_1^\gamma, \pi_2^\gamma \in \Ao nF^\roman{fin}$, and we may say that $\pi_1\equiv \pi_2 \pmod\ell$ if $\pi_1^\gamma \equiv \pi_2^\gamma \pmod \ell$. This relation on complex cuspidal representations of $G$ is then independent of the choice of $\gamma$. 
\endremark 
\head \Smc 
4. Main Theorem 
\endhead 
We prove the Main Theorem of the introduction. As there, $\scr W_F$ denotes the Weil group of a separable algebraic closure $\bar F/F$ of $F$. 
\subhead 
4.1 
\endsubhead 
Let $\Go nF$ denote the set of equivalence classes of irreducible smooth $\Q$-representations of $\scr W_F$ of dimension $n$, and $\Gfin nF$ the set of $\sigma\in \Go nF$ such that $\det\sigma$ has finite order. For $\sigma\in \Gfin nF$, the group $\sigma(\scr W_F)$ is finite and $\sigma$ is effectively a representation of the Galois group of $\bar F/F$ ({\it cf\.} \cite{4} 28.6). We may use the method of 1.1 to attach to $\sigma$ an isomorphism class $[\sigma]_\ell^\roman{ss}$ of smooth, $n$-dimensional semisimple representations of $\scr W_F$. For $\sigma_1,\sigma_2\in \Gfin nF$, we say that $\sigma_1\equiv \sigma_2 \pmod\ell$ if $[\sigma_1]_\ell ^\roman{ss} = [\sigma_2]_\ell^\roman{ss}$. 
\subhead 
4.2 
\endsubhead 
We recall some results from \S1 of \cite{6}. Let $\scr P_F$ be the wild inertia subgroup of $\scr W_F$, and let $\wP F$ denote the set of equivalence classes of irreducible, smooth $\Q$-representations of $\scr P_F$. The group $\scr W_F$ acts on $\wP F$ by conjugation. 
\proclaim{Lemma 1} 
Let $\alpha\in \wP F$. 
\roster 
\item 
The $\scr W_F$-isotropy group of $\alpha$ is of the form $\scr W_E$, for a finite, tamely ramified field extension $E/F$. 
\item 
There exists a smooth representation $\rho$ of $\scr W_E$ such that $\rho|_{\scr P_F} \cong\alpha$. 
\endroster 
\endproclaim 
\demo{Proof} 
Parts (1) and (2) are given by \cite {6} 1.2 Proposition and \cite{6} 1.3 Proposition respectively. \qed 
\enddemo 
In part (1) of Lemma 1, we use the notation $E = Z_F(\alpha)$. 
\par 
Let $\sigma$ be an irreducible, smooth $\Q$-representation of $\scr W_F$. The restriction $\sigma|_{\scr P_F}$ is then a direct sum of representations $\alpha\in \wP F$, any two of which are $\scr W_F$-conjugate. In other words, $\sigma$ determines a unique $\scr W_F$-orbit $r^1_F(\sigma)\in \scr W_F\backslash \wP F$. 
\proclaim{Lemma 2} 
Let $\sigma_1,\sigma_2\in \Gfin nF$. If $\sigma_1\equiv \sigma_2 \pmod \ell$, then $r^1_F(\sigma_1) = r^1_F(\sigma_2)$. 
\endproclaim 
\demo{Proof} 
Since $\scr P_F$ is a pro-$p$ group, the result follows from 1.4 Proposition. \qed 
\enddemo 
If $m\ge1$ is an integer and $\alpha\in \wP F$, we define $\Gg mF\alpha$ to be the set of equivalence classes of irreducible smooth representations of $\scr W_F$ containing $\alpha$ with multiplicity $m$. For example, if $\bk 1$ denotes the trivial representation of $\scr P_F$, then $\Gg mF{\bk 1}$ is the set of equivalence classes of irreducible, $m$-dimensional, {\it tamely ramified\/} representations of $\scr W_F$. In general, we write $\Ggf mF\alpha$ for the set of $\sigma\in \Gg mF\alpha$ such that $\det\sigma$ has finite order. 
\remark{Notation} 
What we here call $\Gg mF\alpha$ is denoted in \cite{6} by $\Gg mF{\scr O}$, where $\scr O = \scr O_F(\alpha)$ is the $\scr W_F$-orbit of $\alpha$. If $E/F$ is a finite separable extension of $F$ with $E\subset \bar F$, we denote by $\Ind_{E/F}$ the smooth induction functor $\Ind_{\scr W_E}^{\scr W_F}$. 
\endremark 
\proclaim{Lemma 3} 
Let $\alpha\in \wP F$, let $E = Z_F(\alpha)$, and let $\rho$ be an irreducible smooth representation of $\scr W_E$ such that $\rho|_{\scr P_F} \cong\alpha$. Let $m\ge1$ be an integer. If $\tau \in \Gg mE{\bk 1}$, the representation $\vS_\rho(\tau) = \Ind_{E/F}\,\rho\otimes\tau$ is irreducible, and the map 
$$
\vS_\rho: \Gg mE{\bk 1} \longrightarrow \Gg mF\alpha 
$$ 
is a bijection. If\/ $\det\rho$ has finite order, then $\vS_\rho\big(\Ggf mE{\bk 1}\big) = \Ggf mF\alpha$. 
\endproclaim 
\demo{Proof} 
The first two assertions combine 1.3 Proposition and 1.4 Theorem of \cite{6}. For the final one, we recall that 
$$
\det\vS_\rho(\tau) = \delta \cdot (\det \rho\otimes \tau)|_{F^\times} = \delta\cdot \det\tau^{\dim\rho}\cdot \det\rho^{\dim\tau}\big|_{F^\times}, 
$$
for a certain tamely ramified character $\delta$ of $F^\times$ of order $\le 2$: if $d_{E/F} = \det \Ind_{E/F}\,1_E$, then $\delta$ is $d_{E/F}^{\dim\rho\otimes\tau}$. The assertion follows. \qed 
\enddemo 
\proclaim{Proposition} 
In the context of Lemma 3, suppose that $\det\rho$ has finite order. If $\tau_1,\tau_2\in \Ggf mE{\bk 1}$, then 
$$
\vS_\rho(\tau_1) \equiv \vS_\rho(\tau_2) \pmod\ell \quad \Longleftrightarrow \quad \tau_1 \equiv \tau_2 \pmod\ell. 
$$ 
\endproclaim 
\demo{Proof} 
Set $\sigma_i  = \vS_\rho(\tau_i)$. Since $\scr P_F$ is a pro-$p$ group, the representation $[\alpha]_\ell^\roman{ss}$ is irreducible. Therefore the representation $\tilde\rho = [\rho]^\roman{ss}_\ell$ is irreducible. If $\mu$ is an irreducible $\F$-representation of $\scr W_E$, trivial on $\scr P_F$, 1.5 Proposition asserts that $\nu_\mu = \Ind_{E/F}\,\tilde\rho\otimes \mu$ is irreducible, and that $\nu_\mu$ determines $\mu$. This implies first the relation $[\sigma_i]^\roman{ss}_\ell = \Ind_{E/F}\, \tilde\rho \otimes [\tau_i]_\ell^\roman{ss}$, $i=1,2$, and then the result. \qed 
\enddemo 
\subhead 
4.3 
\endsubhead 
We follow a parallel path on the other side. For this, it will be convenient to use the language of endo-classes of simple characters: the theory is developed in \cite{1}, but the summary in \cite{6} 2.2 will be adequate. 
\par 
If $\pi \in \Ao nF$, the simple characters contained in $\pi$ lie in a single conjugacy class (3.1), and hence in a single endo-class which we denote $\vt(\pi)$. 
\proclaim{Lemma} 
Let $\pi_1,\pi_2\in \Afin nF$. If $\pi_1\equiv \pi_2 \pmod \ell$, then $\vt(\pi_1) = \vt(\pi_2)$. 
\endproclaim 
\demo{Proof} 
Let $\theta_i$ be an m-simple character occurring in $\pi_i$. The irreducible representation $[\pi_i]_\ell$ contains the $\F$-valued character $[\theta_i]_\ell$ of $H^1_\theta$, obtained by composing $\theta$ with the reduction map $\Z^\times \to \F^\times$ (2.2 Proposition). The isomorphism $[\pi_1]_\ell = [\pi_2]_\ell$ implies that the $[\theta_i]_\ell$ intertwine in $G$ in the usual sense: there exists $g\in G$ such that the characters $[\theta_1]_\ell^g$, $[\theta_2]_\ell$ agree on $(H^1_{\theta_1})^g \cap H^1_{\theta_2}$. Since the values of the $\theta_i$ are $p$-power roots of unity, the characters $\theta_1^g$, $\theta_2$ must agree on $(H^1_{\theta_1})^g \cap H^1_{\theta_2}$. In other words, the $\theta_i$ intertwine in $G$. By the Intertwining Theorem of \cite{7}, they are endo-equivalent and $\vt(\pi_1) = \vt(\pi_2)$ as required. \qed 
\enddemo 
Let $\scr E(F)$ denote the set of endo-equivalence classes of $\Q$-valued simple characters over $F$ and, for $\vT\in \scr E(F)$, let $\deg\vT$ be the degree of $\vT$ ({\it cf\.} \cite{6} (2.2.2)). 
\par 
Let $\vT \in \scr E(F)$ have degree dividing $n$, and set $m = n/\deg\vT$. Let $\Aa mF\vT$ be the set of $\pi\in \Ao nF$ such that $\vt(\pi) = \vT$. Choosing an m-simple character $\theta$ in $G$ of endo-class $\vT$, a representation $\pi\in \Ao nF$ lies in $\Aa mF\vT$ if and only if it contains $\theta$. We attach $\theta$ to a simple stratum $[\frak a,\beta]$, and set $P = F[\beta]$. 
\par 
For the next result, we consider a pair of representations $\pi_1,\pi_2\in \Aaf mF\vT = \Afin nF \cap \Aa mF\vT$. Each $\pi_i$ contains $\theta$, so we let $\vL_i$ be the natural representation of $\bk J_\theta$ on the $\theta$-isotypic subspace of $\pi_i$. As in 3.2, we choose $\kappa\in \scr H(\theta)$ with $\kappa(F^\times)$ finite, and write $\vL_i = \kappa \otimes \lambda_i$, $\lambda_i \in \scr T_0(\theta)$. The representations $[\vL_i]_\ell^\roman{ss}$ and $[\lambda_i]^\roman{ss}_\ell$ are all irreducible, by 3.3 Proposition. 
\proclaim{Proposition} 
The representations $\pi_i$ satisfy $\pi_1\equiv \pi_2 \pmod \ell$ if and only if $[\lambda_1]_\ell^\roman{ss} = [\lambda_2]_\ell^\roman{ss}$. 
\endproclaim 
\demo{Proof} 
We have $[\vL_i]_\ell^\roman{ss} = [\kappa]_\ell^\roman{ss} \otimes [\lambda_i]_\ell^\roman{ss}$. Consequently, if $[\lambda_1]_\ell^\roman{ss} = [\lambda_2]_\ell^\roman{ss}$, then $[\vL_1]_\ell^\roman{ss} \cong [\vL_2]_\ell^\roman{ss}$ and 3.4 Proposition implies $\pi_1 \equiv \pi_2 \pmod\ell$. 
\par 
Conversely, suppose that $\pi_1 \equiv \pi_2 \pmod\ell$. The representation $[\vL_i]^\roman{ss}_\ell$ is the natural representation of $\bk J_\theta$ on the $[\theta]_\ell$-isotypic subspace of $[\pi_i]_\ell$, so $[\vL_1]_\ell^\roman{ss} = [\vL_2]_\ell^\roman{ss}$. The relation $[\vL_i]_\ell^\roman{ss} = [\kappa]_\ell^\roman{ss} \otimes [\lambda_i]_\ell^\roman{ss}$ implies $[\lambda_1]_\ell^\roman{ss} = [\lambda_2]_\ell^\roman{ss}$, by 1.5 Proposition. \qed 
\enddemo 
\subhead 
4.4 
\endsubhead 
We prove the Main Theorem in a special case. It will be notationally convenient to use a base field $E$ and consider representations of $\GL mE$, for some $m\ge1$. Also, if $K/F$ is any finite field extension, we denote by $X_1(K)$ the group of tamely ramified characters of $K^\times$. 
\par
Let $\bk 1$ denote the trivial representation of $\scr P_E$ and let $\bk 0_E\in \scr E(E)$ be the endo-class of trivial simple characters. The Langlands correspondence then induces a bijection $\Gg mE{\bk 1} \to \Aa mE{\bk 0_E}$: this is well known but may be regarded as an instance of \cite{6} 6.4 Corollary. We have to prove: 
\proclaim{Proposition} 
If $\sigma_1,\sigma_2\in \Ggf mE{\bk 1}$, then $\sigma_1\equiv \sigma_2 \pmod\ell$ if and only if $\upr L\sigma_1 \equiv \upr L\sigma_2 \pmod\ell$. 
\endproclaim 
\demo{Proof} 
We recall standard parametrizations of the sets $\Gg mE{\bk 1}$, $\Aa mE{\bk 0_E}$. Let $E_m/E$ be an unramified extension of degree $m$, and set $\vD = \Gal{E_m}E$. Let $X_1(E_m)^{\text{\rm $\vD$-reg}}$ denote the set of $\vD$-regular characters $\chi\in X_1(E_m)$. For $\chi\in X_1(E_m)^{\text{\rm $\vD$-reg}}$, the representation $\sigma_\chi = \Ind_{E_m/E}\,\chi$ is irreducible. It lies in $\Gg mE{\bk 1}$ and depends only on the $\vD$-orbit of $\chi$. Indeed, the map 
$$
\align 
\vD\backslash X_1(E_m)^{\text{\rm $\vD$-reg}} &\longrightarrow \Gg mE{\bk 1}, \\ 
\chi &\longmapsto \sigma_\chi, 
\endalign 
$$ 
is a bijection. Under this map, elements of $\Ggf mE{\bk 1}$ correspond to characters of finite order. 
\par 
On the other hand, let $\theta_0$ denote the trivial character of $U^1_\frak m$, where $\frak m = \M m{\frak o_E}$. Following \cite{3} 2.2, we have a canonical bijection 
$$
\align 
\vD\backslash X_1(E_m)^{\text{\rm $\vD$-reg}} &\longrightarrow \scr T(\theta_0), \\ 
\chi &\longmapsto \gamma_E(\chi), 
\endalign 
$$ 
extending the Green parametrization. The set $\scr T(\theta_0)$ is in canonical bijection with $\Aa mE{\bk 0_E}$ (as in 3.1 Proposition), so we have a further bijection 
$$
\align 
\vD\backslash X_1(E_m)^{\text{\rm $\vD$-reg}} &\longrightarrow \Aa mE{\bk 0_E}, \\ 
\chi &\longmapsto \vG_E(\chi).  
\endalign 
$$ 
\indent 
We examine the behaviour of these constructions relative to congruences. A character $\chi\in X_1(E_m)$, of finite order, may be written uniquely in the form $\chi = \chi_\roman{reg}\,\chi_\roman{sing}$, where $\chi_\roman{reg}$ has order relatively prime to $\ell$ and $\chi_\roman{sing}$ has order a power of $\ell$. 
\proclaim{Lemma} 
Let $\chi,\xi \in X_1(E_m)^{\text{\rm $\vD$-reg}}$ have finite orders. The following conditions are equivalent: 
\roster 
\item $\chi_\roman{reg}$ is $\vD$-conjugate to $\xi_\roman{reg}$; 
\item $\sigma_\chi \equiv \sigma_\xi \pmod\ell$; 
\item $[\gamma_E(\chi)]^\roman{ss}_\ell = [\gamma_E(\xi)]_\ell^\roman{ss}$; 
\item $\vG_E(\chi)\equiv \vG_E(\xi) \pmod\ell$. 
\endroster 
\endproclaim 
\demo{Proof} 
Let $\tilde\chi$ be the character obtained by composing $\chi$ with the reduction map $\Z^\times \to \F^\times$. In particular, $\tilde\chi = (\chi_\roman{reg})\sptilde$, so $\tilde\chi$ determines $\chi_\roman{reg}$ completely. The restriction of $[\sigma_\chi]_\ell^\roman{ss}$ to $\scr W_{E_m}$ is $\sum_{\delta\in \vD} \tilde\chi^\delta$. The same applies to $\sigma_\xi$, so (2) $\Rightarrow$ (1). 
\par 
The decomposition map $*\mapsto [*]_\ell^\roman{ss}$ commutes with induction (\cite{9} \S21B). Identifying $\tilde\chi$ with $[\chi]_\ell^\roman{ss}$, we therefore have $[\sigma_\chi]_\ell^\roman{ss} = \Ind_{E_m/E}\, [\chi]_\ell^\roman{ss}$. That is, (1) $\Rightarrow$ (2). 
\par 
The equivalence of (3) and (4) is a case of 3.4 Proposition, so we consider the relation between (1) and (3). The representations $\gamma_E(\chi)$, $\gamma_E(\xi)$ agree on $E^\times$ if and only if the characters $\chi|_{E^\times}$, $\xi|_{E^\times}$ are equal. The same applies to their reductions modulo $\ell$. Further, $\gamma_E(\chi)$, $\gamma_E(\xi)$ are equivalent on $\GL m{\frak o_E}$ if and only if the restrictions $\chi|_{U_{E_m}}$, $\xi|_{U_{E_m}}$ are $\vD$-conjugate.  Appealing to \cite{17} Theorem 3.6, the reductions modulo $\ell$ of these restrictions are equivalent if and only if $\tilde\chi|_{U_{E_m}}$ is $\vD$-conjugate to $\tilde\xi |_{U_{E_m}}$. Thus (1) is equivalent to (3). \qed 
\enddemo 
To conclude the proof, we need only recall (from, for example, 2.4 Theorem 2 of \cite{5}) that 
$$
\upr L\sigma_\chi = \vG_E(\omega \chi), \quad \chi\in X_1(E_m)^{\text{\rm $\vD$-reg}}, 
$$
where $\omega$ is the unramified character of $E_m^\times$ such that 
$$ 
\omega(x) = (-1)^{(m-1)\ups_{E_m}(x)}, 
$$ 
$\ups_{E_m}$ being the normalized valuation $E_m^\times \to \Bbb Z$. \qed 
\enddemo 
\subhead 
4.5 
\endsubhead 
We return to $\vT\in \scr E(F)$, $\theta$ and $m$, as in 4.3. We choose a simple stratum $[\frak a,\beta]$ in $\M nF$ to which $\theta$ is attached and set $P = F[\beta]$. We let $E/F$ be the maximal tamely ramified sub-extension of $P/F$. We choose $\kappa\in \scr H(\theta)$ and use it to define a bijection $\Aa mE{\bk 0_E} \to \Aa mF\vT$. 
\par 
The lemma of 3.3 yields a bijection 
$$
\vP^P_\kappa: \Aa mP{\bk 0_P} \longrightarrow \Aa mF\vT, 
$$ 
as follows. Let $G_P$ be the $G$-centralizer of $P^\times$ and let $\theta_P$ denote the trivial character of $J^1_\theta\cap G_P$. Thus $\theta_P$ is a trivial m-simple character in $G_P$ and, if $\mu\in \Aa mP{\bk 0_P}$, then $\mu$ contains a unique $\lambda_P\in \scr T(\theta_P)$. In the notation of 3.3, there is a unique representation $\lambda\in \scr T_0(\theta)$ such that $\lambda|_{\bk J_\theta\cap G_P} = \lambda_P$. The tensor product $\lambda\otimes \kappa$ then lies in $\scr T(\theta)$ and we set 
$$ 
\vP^P_\kappa(\mu) = \cind_{\bk J_\theta}^G\,\lambda\otimes\kappa. 
$$
Next, let $P_m/P$ be unramified of degree $m$, and set $\vD = \Gal {P_m}P$. As in 4.4, we have the canonical bijection $\vG_P:\vD\backslash X_1(P_m)^{\text{\rm $\vD$-reg}} \to \Aa mP{\bk 0_P}$. Let $E_m/F$ be the maximal tamely ramified sub-extension of $P_m/F$, so that $E_m/E$ is unramified of degree $m$. We identify $\vD$ with $\Gal {E_m}E$ by restriction. The extension $P_m/E_m$ is totally wildly ramified, so composition with the field norm $\N{P_m}{E_m}$ induces a $\vD$-isomorphism $X_1(E_m) \to X_1(P_m)$. This in turn induces a bijection 
$$ 
\beta_{P/E}: \vD\backslash X_1(E_m)^{\text{\rm $\vD$-reg}} \longrightarrow \vD\backslash X_1(P_m)^{\text{\rm $\vD$-reg}}. 
$$ 
We so obtain a bijection 
$$
\align 
\roman b_{P/E}: \Aa mE{\bk 0_E} &\longrightarrow \Aa mP{\bk 0_P}, \\ 
\mu &\longmapsto \vG_P\circ\beta_{P/E}\circ\vG_E^{-1}(\mu). 
\endalign 
$$ 
We now define 
$$
\align 
\vP_\kappa: \Aa mE{\bk 0_E} &\longrightarrow \Aa mF\vT, \\ 
\mu&\longmapsto \vP_\kappa^P(\roman b_{P/E}\,\mu). 
\endalign 
$$ 
The map $\vP_\kappa$ is indeed a bijection and, if we choose $\kappa$ so that $\kappa(F^\times)$ is finite, then $\vP_\kappa$ maps $\Aaf mE{\bk 0_E}$ onto $\Aaf mF\vT$. 
\proclaim{Proposition} 
Let $\kappa\in \scr H(\theta)$, and suppose that $\kappa(F^\times)$ is finite. If $\mu_1,\mu_2 \in \Aaf mE{\bk 0_E}$, then 
$$
\mu_1\equiv\mu_2 \pmod\ell \quad \Longleftrightarrow \quad \vP_\kappa(\mu_1) \equiv \vP_\kappa(\mu_2) \pmod\ell. 
$$
\endproclaim 
\demo{Proof} 
The proposition is just a re-statement of 4.3 Proposition. \qed 
\enddemo 
\subhead 
4.6 
\endsubhead 
We prove the Main Theorem in the general case. Let $\sigma_1,\sigma_2\in \Gfin nF$, and set $\pi_i = \upr L\sigma_i\in \Ao nF^\roman{fin}$, $i=1,2$. 
\par 
We dispose of some trivial cases. If $r_F^1(\sigma_1) \neq r^1_F(\sigma_2)$, then $\sigma_1\not\equiv \sigma_2 \pmod \ell$, by 4.2 Lemma 2. This implies, via 6.1 Theorem of \cite{6}, that $\vt(\pi_1)\neq \vt(\pi_2)$, whence $\pi_1 \not\equiv \pi_2 \pmod \ell$ (4.3 Lemma). This argument is reversible, so we need only consider the case $r^1_F(\sigma_1) = r^1_F(\sigma_2)$. In other words, there exists $\alpha\in \wP F$ such that $\sigma_1,\sigma_2\in \Ggf mF\alpha$, for some integer $m\ge1$. In light of 4.4, we need only treat the case where $\alpha$ is not trivial. We set $E = Z_F(\alpha)$. 
\par
According to 6.4 Corollary of \cite{6}, the set $\upr L(\Gg mF\alpha)$ is of the form $\Aa mF\vT$, for a unique $\vT\in \scr E(F)$. In other words, $\pi_i\in \Aaf mF\vT$, $i=1,2$. As in 4.5, we choose an m-simple character $\theta$ in $G = \GL nF$ of endo-class $\vT$, and attach it to a simple stratum $[\frak a,\beta]$ with $P = F[\beta]$. According to the Tame Parameter Theorem of \cite{6} 6.3, the maximal tamely ramified sub-extension of $P/F$ is isomorphic to $E/F$. We fix such an isomorphism, and henceforward view $E$ as a subfield of $P$. We choose a smooth representation $\rho$ of $\scr W_E$ such that $\rho|_{\scr P_F} \cong \alpha$ (as we may, by 4.2 Lemma 1) and also a representation $\kappa\in \scr H(\theta)$. We have a (non-commutative) diagram of bijections,   
$$
\CD 
\Gg mE{\bk 1} @>{\vS_\rho}>> \Gg mF\alpha \\ 
@VVV @VVV \\ 
\Aa mE{\bk 0_E} @>>{\vP_\kappa}> \Aa mF\vT 
\endCD \tag 4.6.1 
$$ 
in which the vertical arrows are given by the Langlands correspondence. 
\proclaim{Lemma} 
Let $\rho$ be a smooth representation of $\scr W_E$ such that $\rho|_{\scr P_F} \cong \alpha$. There exists a unique $\kappa = \kappa(\rho)\in \scr H(\theta)$ such that the diagram \rom{(4.6.1)} commutes. If $\det\rho$ has finite order, then $\kappa(F^\times)$ is finite. 
\endproclaim 
\demo{Proof} 
The first assertion is a restatement of the Comparison Theorem of \cite{6}, while the second is immediate. \qed 
\enddemo 
We take $\rho$, $\kappa$ as in the lemma, with $\det\rho$ of finite order. For $\sigma_i \in \Ggf mF\alpha$, we write $\sigma_i = \vS_\rho(\tau_i)$, $\tau_i\in \Gg mE{\bk 1}$, $i=1,2$. The condition $\sigma_1\equiv \sigma_2 \pmod \ell$ is equivalent to $\tau_1 \equiv \tau_2 \pmod\ell$, by 4.2 Proposition. This, in turn, is equivalent to $\upr L\tau_1 \equiv \upr L\tau_2 \pmod \ell$, by 4.4 Proposition. The lemma informs us that $\pi_i = \vP_\kappa(\upr L\tau_i)$, $i=1,2$, while 4.5 Proposition asserts that $\pi_1\equiv \pi_2 \pmod\ell$ if and only if $\upr L\tau_1 \equiv \upr L\tau_2 \pmod \ell$. This completes the proof of the Main Theorem. \qed 
 
\enddocument